\documentclass[10pt]{article}

\usepackage{latexsym}
\usepackage{amsmath, amsthm}
\usepackage{amsfonts}
\usepackage{bbm}
\usepackage{pdfsync}
\usepackage{graphicx}

\newtheorem{theorem}{Theorem}
\newtheorem{prop}{Proposition}
\newtheorem{lemma}{Lemma}

\numberwithin{equation}{section}
\theoremstyle{plain}

\def\Bbb#1{{\mathbb#1}}
\def\RR{\Bbb{R}}
\def\NN{\Bbb{N}}

\def\EE{\Bbb{E}}
\def\PP{\Bbb{P}}

\def\B{{\cal B}}
\def\x{{\mathbf{x}}}
\def\X{{\mathbf{X}}}

\def\eps{\varepsilon}

\def\half{\frac{1}{2}}

\begin{document}

\begin{center}
\vspace*{5mm}

{\bf \Large Conditional Distribution of Heavy Tailed Random Variables on Large Deviations of their Sum.}\\
\vspace{6mm}
In\'{e}s Armend\'{a}riz\footnote[1]{Universidad de San Andr\'{e}s, Vito Dumas 284, B1644BID, Victoria,  Argentina. E-mail:{\tt iarmendariz@udesa.edu.ar}},\
Michail Loulakis\footnote[2]{Department of Applied Mathematics, University of Crete, and Institute of Applied and Computational Mathematics, FORTH,  Crete.\ Knossos Avenue, 714 09 Heraklion Crete, Greece. E-mail: {\tt loulakis@tem.uoc.gr}}
\end{center}
\vspace{6mm}
\small
\hspace{1cm} \begin{minipage}[t]{4.9in}
ABSTRACT: It is known that large deviations of sums of subexponential random variables are most likely realised by deviations of a single random variable. In this article we give a detailed picture of how subexponential random variables are distributed when a large deviation of the sum is observed.\\
\\
{\em AMS 2000 Mathematics Subject Classification}: 60F10 \\
\\
{\em Keywords}: Large Deviations, Subexponential Distributions, Conditional Limit Theorem, Gibbs conditioning principle.
\end{minipage}
\normalsize
\vspace{1mm}
\section{Introduction}
Let $X_1,X_2,\ldots$ be a sequence of i.i.d. random variables with common distribution $\mu$ defined on a probability space $(\Omega,{\cal F},\PP)$, and let
\[
S_n=X_1+\cdots+X_n, \qquad n\ge 1.
\]
The most classical problem in large deviations is establishing asymptotic expressions for
\begin{equation}
\bar{F}_n(x):=\PP\big[ S_n > x \big]
\label{tails}
\end{equation}
when this quantity converges to zero. The answer depends heavily on the nature of the tails of the distribution $\mu$. When the moment generating function is finite in a neighborhood of the origin (Cram\'{e}r's condition) Cram\'{e}r derived asymptotic expressions for $\bar{F}_n(x)$ valid uniformly over different ranges of $x$-values. These results were later refined by Petrov (cf. \cite{MR1312369}.) In this case Gibbs conditioning principle provides an answer to how a large deviation of the sum is typically realised: subject to the large deviation, the random variables $\{X_i\}$ become independent in the limit, but their marginal distribution is modified in such a way that the behavior imposed on the sum now becomes typical. In particular, no single random variable becomes excessively large compared to the others.\\
\\
The situation is totally different when Cram\'{e}r's condition is violated. It is known since the classical works of Heyde \cite{MR0240854} and Nagaev \cite{MR0282396} that large deviations of sums of independent heavy-tailed random variables are typically realised by one random variable taking a very large value.  In this article we investigate the conditional distribution of the random variables $\{X_i\}_{1\le i\le n}$ subject to a large deviation of their sum $S_n$.  It turns out that as $n\to \infty$ this conditional distribution converges to a product of $n-1$ copies of $\mu$, while the remaining variable realises the large deviation event by taking a very large value. We determine when the fluctuations around that value have a scaling limit, and we show that given the sum exceeds a large value, the maximum is asymptotically independent of the smallest variables. \\
\\
\section{Notation and Results}
Let $\{X_n\}_{n\in\NN}$ be a sequence of i.i.d. random variables with common distribution $\mu$ defined on a probability space $(\Omega,{\cal F},\PP)$.  We denote by $F$ their distribution function $F(x)=\mu(-\infty,x].$ We are interested in the case where $F$ is in the class of subexponential distributions, that is
\begin{equation}
\lim_{x\to\infty}\frac{\bar{F}(x+y)}{\bar{F}(x)}=1,\qquad \forall y\in\RR,
\label{longtail}
\end{equation}
with $\bar{F}(x)=\PP\big[X_k>x\big]$, and
\begin{equation}
\lim_{x\to\infty}\frac{\bar{F}_n(x)}{n\bar{F}(x)}=1,\qquad \forall n\in\NN,
\label{subexp}
\end{equation}
where $\bar{F}_n(x)$ is defined in (\ref{tails}). If $\mu$ is supported on the positive half-line then (\ref{longtail}) is implied by (\ref{subexp}), and in that case subexponentiality can be defined by the latter condition alone.  Since it is generally true that for all $n\in\NN$ we have
\[
\lim_{x\to\infty}\frac{\PP\big[\max_{1\le k\le n} X_k>x\big]}{n\bar{F}(x)}=1,
\]
equation (\ref{subexp}) states that the tail of the sum in a sample of independent $\mu$-distributed random variables is determined by the tail of the largest variable. These distributions arise naturally when modelling heavy-tailed phenomena. For instance, individual claims in insurance or large interarrival times in queuing systems are usually modelled by distributions of this kind. Typical members of this class include distributions with regularly varying, lognormal-type, or Weibull-type tails. Sufficient conditions for a given distribution to be subexponential that are straightforward to check can be found in \cite{MR0391222}.\\
\\
An immediate consequence of (\ref{subexp}) is the existence of a sequence $d_n\to\infty$ such that
\begin{equation}
\lim_{n\to\infty}\sup_{x\ge d_n}\left|\frac{\bar{F}_n(x)}{n\bar{F}(x)}-1\right|=0.
\label{asymptotic}
\end{equation}
\noindent
A large amount of work has been done for determining a threshold $d_n$, for which (\ref{asymptotic}) holds. Interested readers can find reviews on the topic in \cite{MR542129,MR1312369}. A very nice account is also provided by Mikosch and A. Nagaev in \cite{MR1652936}. Denisov, Dieker and Schneer give an up-to-date treatment of this problem in  \cite{MR2440928}.\\ 
\\
\noindent
When the distribution $\mu$ satisfies a local version of (\ref{subexp}), a local version of (\ref{asymptotic}) is valid. Let $\Delta=(0,s]$ for some $s>0$ and denote by $x+\Delta$ the interval $(x,x+s]$. We say that $\mu$ is $\Delta$-subexponential if $\mu\big[x+\Delta\big]>0$ for all sufficiently large $x$ and
\begin{equation}
\lim_{x\to\infty}\frac{\mu\big[x+y+\Delta\big]}{\mu\big[x+\Delta\big]}=1,\qquad \forall y\in\RR,
\label{locallongtail}
\end{equation}
and
\begin{equation}
\lim_{x\to\infty}\frac{\PP\big[S_n\in x+\Delta\big]}{n\mu\big[x+\Delta\big]}=1,\qquad \forall n\in\NN.
\label{localsubexp}
\end{equation}
\noindent
The concept of $\Delta$-subexponentiality was introduced in \cite{MR1982040} by Asmussen, Foss and Korshunov. A $\Delta$-subexponential distribution is also $m\Delta$-subexponential for all $m\in\mathbb{N}$, $m\Delta=(0,ms]$, and subexponential in the sense of (\ref{longtail}) and (\ref{subexp}) (cf. \cite{MR1982040}.) Even though there are examples of subexponential distributions that are not $\Delta$-subexponential for finite $\Delta$, most distributions  that are used in practice are. This can be easily verified using the sufficient conditions for $\Delta$-subexponentiality provided in \cite{MR1982040}. The asymptotics for the large deviation probabilities are now given by
\begin{equation}
\lim_{n\to\infty}\sup_{x\ge d_n}\left|\frac{\PP\big[S_n\in x+\Delta\big]}{n\mu\big[x+\Delta\big]}-1\right|=0,
\label{localasymptotic}
\end{equation}
\noindent
and sufficient conditions on $d_n$ for (\ref{localasymptotic}) to hold can be found in \cite{MR2440928}. \\
\\
We are interested in the conditional distribution of the variables $\{X_n\}$ subject to a large deviation of their sum. Assuming that (\ref{localasymptotic}) holds for some interval $\Delta$ that may be finite or infinite, we would like to determine the asymptotic behaviour of
\begin{align*}
&\mu_{n,x}^{\Delta}\big[A\big]=\PP\big[(X_1,\ldots,X_n)\in A\ \big|\  S_n\in x+\Delta\big]
\end{align*}
when $n\to\infty$ and $x\ge d_n$.  Note that when $\Delta=(0,\infty)$ the definition of $\Delta$-subexponentiality reduces to the standard definition of subexponentiality and (\ref{localasymptotic}) reduces to (\ref{asymptotic}). This allows to treat rare events of the form $\{S_n\in(x,x+s]\, \}$ or $\{S_n>x\}$ simultaneously in Theorem \ref{thmgen} below.\\
\\
\noindent
A related question was raised in \cite{MR2013129} where certain subexponential families $\mu$ of lattice type are considered under  $\{S_n=x(n)\}$, and it is shown that the finite dimensional marginals of the conditional distribution converge to a product of copies of $\mu$. \\
\\
\noindent
We will denote by $T:\cup_{n\in \NN} \RR^n \to \cup_{n\in \NN} \RR^n$ the operator that exchanges the last and the maximum component of a finite sequence:
\[
T(x_1,\dots,x_n)_k=\begin{cases} \max_{1 \le i\le n} x_i & \text{if } k=n,\\
                           x_n & \text{if } x_k>\max_{1\le i< k} x_i \text{ and } x_k=\max_{i\ge k}x_i,\\
                           x_k & \text{otherwise.}
            \end{cases}
\]
\begin{theorem}
Suppose $\mu$ is $\Delta$-subexponential. There exists a sequence $q_n$ such that
\[
\lim_{n\to\infty}\sup_{x\ge\, q_n}\sup_{A\in\B(\RR^{n-1})}\left|\mu_{n,x}^{\Delta}\circ T^{-1}\big[A\times\RR \big]-\mu^{n-1}\big[A\big]\right|=0.
\]
\label{thmgen}
\end{theorem}
\noindent
The sequence $q_n$ in the statement can be easily computed from the sequence $d_n$ in (\ref{localasymptotic}) and $F$, and in most cases turns out to be $d_n$ itself. Note that $\mu_{n,x}^{\Delta}\circ T^{-1}\big[A\times\RR \big]$ is the measure assigned to $A\in\RR^{n-1}$ by the conditional distribution of the $n-1$ smallest variables. In other words, Theorem \ref{thmgen} states that under (\ref{asymptotic}), conditioning on $\{S_n\in x+\Delta\}$ affects only the maximum in the limit, and the $n-1$ smallest variables become asymptotically independent. Such a result is rather uncommon, and when $\mu$ satisfies Cram\'er's condition an analogous statement is not true. Now, any limit theorem for  i.i.d. random variables with distribution $\mu$ can be cast in this setting. For instance, we could obtain conditional limit theorems for the statistics of any order $k>1$: the $k$-th order statistic of $(X_1,\ldots,X_n)$ subject to the condition $S_n\in x+\Delta,\, x\ge q_n$ asymptotically behaves like the $(k-1)$-th order statistic of an independent sample.\\
\\
Unlike the asymptotic independence of the smallest variables, the
fluctuations of the maximum $M_n=\max_{1\le i\le n} X_i$ 
and its dependence on the smallest variables are influenced by the
form of conditioning. When $\Delta=(0,s]$ the condition we impose on
the sum is very restrictive and the fluctuations of the maximum are
determined by the fluctuations of the sum of the smallest variables.
This can be easily seen since
$\mu_{n,x}^{\Delta}\big[M_n+\sum_{j=1}^{n-1}(TX)_i\in (x,x+s]\,
\big]=1$ by definition. Therefore, if the (unconditioned)
distribution of $S_{n-1}/b_{n}$ converges to a stable law $H$, it follows immediately
from Theorem 1 that under $\mu_{n,x}^{\Delta}$ we have
\begin{equation}
\frac{M_n-x}{b_{n}} \stackrel{d}{\longrightarrow} -H.
\label{localmaxfluct}
\end{equation}
Note that the converse is also true. In particular, the fluctuations of the conditional maximum are typically two-sided and they have a non trivial scaling limit if and only if $\mu$ is attracted to a stable distribution. In \cite{MR2520125} Theorem 1 is proved for a particular family of lattice distributions subject to $\{S_n=x\}$ and this observation is used to obtain a limit theorem for the fluctuations of the maximum in a system of interacting particles.\\
\\
On the other hand when we condition on $\{S_n>x\}$ it turns out that the maximum coordinate is asymptotically independent of the smallest variables, its fluctuations around $x$ are one-sided, and they have a non trivial scaling limit if and only if $\mu$ is in the maximum domain of attraction of an extreme value distribution. For ease of notation, we will now drop $\Delta=(0,\infty]$ from the notation,
\[
\mu_{n,x}\big[A\big]=\PP\big[(X_1,\ldots,X_n)\in A\ \big|\  S_n> x\big].
\]
Let $\nu_{x}$ stand for the conditional distribution of $X_i$ subject to $X_i>x$. That is,
\[
\nu_x\big[A\big]=\PP\big[X_i\in A\, \big|\, X_i>x\big]= \frac{\mu\big[A\cap(x,\infty)\big]}{\bar{F}(x)}.
\]
We will use $\|\nu\|_{t.v.}$ to denote the total variation norm of a signed Borel measure on $\RR^n$. That is
\[
\|\nu\|_{t.v.}=\sup_{A\in{\cal B}(\RR^n)}\big|\nu\big(A\big)\big|.
\]
\begin{theorem} Suppose $\mu$ is subexponential. Then
\begin{equation}
\lim_{n\to\infty}\sup_{x\ge q_n} \left\|\mu_{n,x}\circ T^{-1}- (\mu^{n-1}\times \nu_x)\right\|_{t.v.}=0,
\label{thm2assert}
\end{equation}
where $q_n$ is the sequence appearing in Theorem \ref{thmgen}.
\label{thmrv}
\end{theorem}
\noindent
Since the distribution of $(X_1,\dots,X_n)$ subject to $\{S_n>x\}$ is
clearly exchangeable, the position of the maximum coordinate is
uniformly distributed among $1,\dots,n$. Theorem \ref{thmrv}
states that the conditional distribution of the maximum coordinate becomes asymptotically
a randomly located $\nu_x$, while the law of the remaining
$n-1$ variables is the product $\mu^{n-1}$ as was established in Theorem
\ref{thmgen}.\\
\\
It is interesting to examine whether (\ref{thm2assert}) entails a limit theorem for the fluctuations of the maximum around $x$, that is, whether there exists a scaling function $\psi(\cdot)$ such that under $\mu_{n,x}$ we have
\begin{equation}
\frac{M_n-x}{\psi(x)} \stackrel{d}{\longrightarrow} \Lambda,
\label{maxfluct}
\end{equation}
for some non trivial distribution $\Lambda$. In view of Theorem \ref{thmrv} this is equivalent to asking when
\begin{equation}
\nu_x\big[(x+u\psi(x),\infty)\big]=\frac{\bar{F}\big((x+u\psi(x)\big)}{\bar{F}(x)}
\label{residual}
\end{equation}
converges as $x\to\infty$ to a nontrivial function of $u$. This is
precisely the subject of \cite{MR0359049}, where Balkema and de Haan
determine all possible scaling limits of residual life times as the
survival time goes to infinity, and the corresponding domains of
attraction. It follows from their results (Theorems 1, 3 and 4
there) that nontrivial limits in the right hand side of
(\ref{maxfluct}) can only be of two types.
\begin{enumerate}
\item
An exponential distribution of rate 1 if and only if $\mu$ is in the maximum domain of attraction of the Gumbel distribution. In this case $\psi$ can be determined by requiring the expression in (\ref{residual}) to converge to $e^{-u}$.
\item
A Pareto destribution on $\mathbb{R}_{+}$ with $\bar{\Lambda}(u)=(1+u)^{-\alpha}$ and $\alpha>0$, if and only if $\mu$ has regularly varying tails with index $-\alpha$, that is $\bar{F}(x)=x^{-\alpha}L(x)$, as $\ x\to\infty$, and $L$ is a slowly varying function. Note that this is equivalent to $\mu$ being in the maximum domain of attraction of the Fr\'echet distribution with index $\alpha$ (cf. \cite{MR1015093}.) In this case $\psi(x)=x$.
\end{enumerate}
\noindent
That regularly varying distributions satisfy our main assumption (\ref{asymptotic}) is a long known fact (cf. \cite{MR0240854,MR0282396}.) In particular, if $\alpha>2$ one can choose $d_n=\sqrt{tn\log n}$ for any $t>\alpha-2$  (cf. \cite{MR542129}.)  The articles \cite{MR2440928,MR1652936} are excellent references for subexponential distributions in the maximum domain of attraction of the Gumbel distribution and the corresponding sequences $d_n$ for which (\ref{asymptotic}) holds.\\
\\
\noindent
{\bf Remarks}\\
1. Theorem \ref{thmrv} and the discussion following it generalise a result of Mikosch and Nagaev (Proposition 4.4 in \cite{MR1652936}) where they prove (\ref{maxfluct}) for some of the most commonly used subexponential distributions. \\
\\
\noindent
2. Theorem \ref{thmrv} also generalises an old result by Richard Durrett \cite{MR576888}. In that article it is proved that if $\mu$ is regularly varying with index $\alpha<-2$ and $\EE\big[X_1\big]=-b<0$,  then
\[
\big(\ \frac{S_{[n\cdot]}}{n}\ \big|\ S_n>0\ \big) \Rightarrow J_{\alpha,b}1_{\{U\le \cdot\}}-b\ \cdot,
\]
where $U$ is uniform in [0,1] and $J_{\alpha,b}$ is independent of $U$ with Pareto distribution. 
As in Corollary 3 in \cite{MR2520125}, Theorem \ref{thmrv} also establishes a
conditional invariance principle for the sum of the random variables
cut-off at a level $\eps n$. \\
\\
\noindent
Simple modifications of the proofs give ramifications of Theorems  \ref{thmgen} and \ref{thmrv}. For instance, the results remain true if $n$ is fixed and we only let $x\to\infty$. Theorem \ref{thmgen} remains valid if we let the size of $\Delta$ grow to infinity with $n$ and $x$, and a variant of Theorem \ref{thmrv} is satisfied if $\Delta$ grows fast enough. This enables to explore precisely how the fluctuations of the maximum switch from a stable to a residual life time nature. These ramifications are discussed in section 4, after the proof of the Theorems that follows.

\section{Proof of the Theorems}
\label{Proofs}

\noindent
Given a vector $\mathbf{x}=(x_1,\ldots,x_n)\in\RR^n$ we denote by
$M_{\x}$  the coordinate of maximum size and by $m_{\mathbf{x}}$  its position. Precisely,
\[
M_{\x}=\max_{1\le k \le n} x_i
\qquad \mbox{and}\qquad m_{\mathbf{x}}=k \Leftrightarrow x_k> x_j, \ j< k,\,\,\mbox{and}\,\,\, x_k\ge x_j,\, j\ge k.
\]
\noindent
We will denote by $\sigma^j$ the operator that exchanges the $j$-th and the last coordinate of $\mathbf{x}$, that is
\[
\sigma^j(x_1,\ldots,x_j,\ldots,x_{n-1},x_n)=(x_1,\ldots,x_n,\ldots,x_{n-1},x_j).
\]
\noindent
With this notation we may write $T\mathbf{x}=\sigma^{m_{\mathbf{x}}}\mathbf{x}$. Let also $\X=(X_1,\dots,X_n)$, $\X^{n-1}=(X_1,\dots,X_{n-1})$.
We begin with some elementary observations that will be useful for both proofs.\\
\\
The convergence in (\ref{locallongtail}) is in fact uniform over compact $y$--sets. This follows from the uniform convergence theorem for
slowly varying functions (see \cite{MR1015093}, Theorem 1.2.1), as (\ref{locallongtail}) implies that $x\mapsto \mu\big[\log x+\Delta\big]$ is slowly varying. In particular, if $b_n$ is any sequence growing to infinity there exists a sequence $m_n\to\infty$ such that
\begin{equation}\label{uct}
\lim_{n\to\infty}\sup_{x\ge m_n}\sup_{0\le y\le b_n}\left|\frac{\mu\big[x-y+\Delta\big]}{\mu\big[x+\Delta\big]}-1\right|=0.
\end{equation}
This in turn implies that there exists a sequence $\ell_n\gg b_n$ such that
\begin{align}\label{ratio}
D_n(L):=\sup_{x\ge \ell_n}\sup_{|y|\le Lb_n}\left(1-\frac{\mu\big[x-y+\Delta\big]}{\mu\big[x+\Delta\big]}\right)\longrightarrow 0, \qquad\text{as }{n\to\infty}, \qquad\forall L>0.
\end{align}
To see this, iterate (\ref{uct}) using the fact the limit is uniform in $x\ge m_n$ to get
\begin{equation}
\lim_{n\to\infty}\sup_{x\ge m_n}\sup_{-Lb_n\le y\le b_n}\left|\frac{\mu\big[x-y+\Delta\big]}{\mu\big[x+\Delta\big]}-1\right|=0.
\label{arXiversion}
\end{equation}
Now, if $\rho_n$ is any sequence increasing to infinity we may choose $\ell_n=m_n+\rho_n b_n$.\\
\\
The sequence $q_n$ in the statement of the theorems can be constructed as follows. Take a sequence $b_n$ such that $S_{n-1}/b_{n}$ is tight and choose $\ell_n$ so that (\ref{ratio}) holds. We may then choose $q_n=d_n\vee \ell_{n}$. Very often, in fact in all cases we are aware of where a threshold $d_n$ in (\ref{asymptotic}) or (\ref{localasymptotic}) is explicitly known, and certainly for the $d_n$ constructed in \cite{MR2440928}, we can choose $\ell_{n}\le d_n$ so the supremum in the theorems is taken for $x\ge d_n$.
Finally, note that for any $B\in\B(\RR^{n})$ we have
\begin{align}
\PP\big[\,T\mathbf{X}\in B, S_n\in x+\Delta \big] &=\sum_{j=1}^n\PP\big[\,T\mathbf{X}\in B, S_n\in x+\Delta,\ m_\mathbf{X}=j\big] \nonumber\\
&\ge \sum_{j=1}^n\PP\big[\sigma^j\mathbf{X}\in B, S_n\in x+\Delta,\ m_{\sigma^j\mathbf{X}}=n\big]\nonumber\\
&=n\,\PP[\,\mathbf{X}\in B, S_n\in x+\Delta,\ m_\mathbf{X}=n\big].
\label{note1}
\end{align}
The last equality holds because $\PP$ is invariant under $\sigma^j$. The penultimate inequality holds because $m_{\sigma^j\mathbf{x}}=n\Rightarrow m_\mathbf{x}=j$ (notice however that  if $\mu$ is atomless  this inequality and (\ref{note2}) below are in fact equalities.) In view of (\ref{note1}) we have
\begin{equation}
\mu_{n,x}^{\Delta}\circ T^{-1}\big[B\big]\ge\frac{n\,\PP[\,\mathbf{X}\in B,\, S_n\in x+\Delta,\, m_\mathbf{X}=n\big]}{\PP\big[S_n\in x+\Delta\big]}.
\label{note2}
\end{equation}
{\em Proof of Theorem \ref{thmgen}. }
\noindent Consider $A\in \B(\RR^{n-1})$ as in the statement of the
theorem and fix $L\in \NN$. We have
\begin{align}
 &\PP\big[\, \X\in A\times \RR,\,S_n\in x+\Delta,\,m_{\X}=n\big]\ge\,\PP\big[\,
\X\in A\times \RR,\,S_n\in x+\Delta,\, |S_{n-1}|<Lb_{n},\,m_{\X}=n\big]\notag \\
&\qquad\qquad\qquad\ge \PP\big[\, \X\in A \times \RR,\,S_n\in x+\Delta,\, |S_{n-1}|<L
b_{n},\,M_{\X^{n-1}}\le x-Lb_{n}\big]\notag \\
&\qquad\qquad\qquad=\int_{\X^{n-1}\in A\cap G} \,\mu\big[x-S_{n-1}+\Delta\big]\, d\PP,\notag
\end{align}
where
\begin{equation}
G=G(n,L,x)=\big\{\mathbf{u}\in\RR^{n-1}:\, \big|\sum_{i=1}^{n-1} u_i\big|<Lb_n,\, M_{\mathbf{u}}\le x-Lb_n\big\}.
\label{G}
\end{equation}
Notice that when $\mathbf{u}\in G$ and $x\ge \ell_n$ we have
\[
\mu\big[x-\sum_{i=1}^{n-1}u_i+\Delta\big]\ge \big(1-D_n(L)\big)\mu\big[x+\Delta\big],
\]
so that (\ref{note2}) can be reinforced to
\begin{align*}
\mu_{n,x}^{\Delta}\circ T^{-1}\big[A\times\mathbb{R}\big]&\ge\big(1-D_n(L)\big)\,\frac{n\mu\big[x+\Delta\big]}{\PP\big[S_n\in x+\Delta\big]} \, \PP\big[\X^{n-1}\in A\cap G\big]
\end{align*}
\noindent
giving the estimate
\begin{align*}
\mu_{n,x}^{\Delta}\circ T^{-1}\big[\, A\times \RR\big]-\PP\big[\X^{n-1} \in A\big]\,\ge\,-\left( \PP\big[\X^{n-1}\notin G\big]+D_n(L)+\left|\frac{n\mu\big[x+\Delta\big]}{\PP\big[S_n\in x+\Delta\big]}-1\right|\,\right).
\end{align*}
Denote the expression in the parenthesis on the right hand side above by $R(n,L,x)$. We can get an upper bound by applying the same estimate for $\RR^{n-1}\setminus A$, the complement of $A$. Combining the two bounds we get
\[
\left|\, \mu_{n,x}^{\Delta}\circ T^{-1}\big[\, A\times \RR\big]-\PP\big[\, \X^{n-1} \in A\big] \,\right|\le R(n,L,x).
\]
Now the sequence $S_{n-1}/b_n$ is tight, so we have
\begin{align}\label{rc}
\lim_{L\to \infty} \sup_{n}
\PP\big[\,|S_{n-1}|\ge Lb_n\big]=0.
\end{align}
On the other hand, it is known (see \cite{MR0270403}, Section lX.7) that
\begin{equation}
\lim_{L\to \infty} \sup_n \,n[F(-Lb_n)+\bar{F}(Lb_n)]=0,
\label{tighty}
\end{equation}
and since $\ell_n\gg b_n$ we have
\begin{align*}
\sup_{x\ge \ell_n} \PP\big[\,M_{\X^{n-1}} > x-Lb_n\big]=1-\inf_{x\ge \ell_n}\big(1-\bar{F}(x-Lb_n)\big)^{n-1}
\longrightarrow 0, \qquad\text{as } n\to\infty.
\end{align*}
Combining this limit with (\ref{rc}) we see that $\PP\big[\X^{n-1}\notin G\big]$ vanishes uniformly on $x\ge \ell_n$, as $n\to\infty$ and then $L\to\infty$.
The result now follows from (\ref{localasymptotic}) and (\ref{ratio}).
$\Box$\\

\noindent
{\em Proof of Theorem \ref{thmrv}.}
The proof follows the general outline of Theorem \ref{thmgen}.\\
\\
It is sufficient to show that
\begin{equation}
\lim_{n\to\infty}\sup_{x\ge d_n\vee\,\ell_n}\sup_{U\in{\cal R}} \Big|\mu_{n,x}\circ T^{-1}\big[U]-\mu^{n-1}\times\nu_x\big[U\big]\Big|=0,
\label{Teq}
\end{equation}
where the supremum above is taken over the class ${\cal R}$ of finite disjoint unions of rectangles $A_j\times B_j$ with $A_j\in{\cal B}(\RR^{n-1})$ and $B_j\in{\cal B}(\RR)$. Recall the definition of  $G\subset\RR^{n-1}$ in (\ref{G}) and define $I=(x+Lb_n,\infty)$. Using (\ref{note2}) we have
\begin{align*}
\mu_{n,x}\circ T^{-1}\Big[\bigcup_j A_j\times B_j\Big] &=\sum_j \mu_{n,x}\circ T^{-1}\big[A_j\times B_j\big]\\
&\ge \sum_j \frac{n\PP\big[\X^{n-1}\in A_j,\, X_n\in B_j,\, S_n>x,\, m_{\X}=n\big]}{\bar{F}_n(x)}\\
&\ge\frac{n}{\bar{F}_n(x)}\sum_j\PP\big[\X^{n-1}\in A_j\cap G,\, X_n\in B_j\cap I\big]\\
&=\frac{n\bar{F}(x)}{\bar{F}_n(x)}\mu^{n-1}\times\nu_x\Big[\big(\bigcup_j A_j\times B_j\big)\cap (G\times I)\,\Big].
\end{align*}
Just as in the proof of Theorem \ref{thmgen} this gives that for all $U\in{\cal R}$ we have
\begin{align*}
\hspace{-2cm}\mu_{n,x}\circ T^{-1}\big[U\big] -\mu^{n-1}\times\nu_x\big[U\big]&\ge -\Big(\mu^{n-1}\times\nu_x\Big[(G\times I)'\Big]+\big|\frac{n\bar{F}(x)}{\bar{F}_n(x)}-1\big|\,\Big)\\
&\ge -\Big(\mu^{n-1}\big[G'\big]+\nu_x\big[I'\big]+\big|\frac{n\bar{F}(x)}{\bar{F}_n(x)}-1\big|\,\Big)\\
&\ge -R(n,L,x),\qquad \forall x\ge \ell_n.
\end{align*}
In the previous equation and in the following the prime symbol denotes the complement of a set in the appropriate space: $(G\times I)'=\RR^n \setminus (G\times I)$, $I'=\RR \setminus I$ and
$G'=\RR^{n-1} \setminus G$.\\

\noindent
Since ${\cal R}$ is closed under complementation we can also get an upper bound by applying the previous inequality for $U'$ to get
\[
\big|\mu_{n,x}\circ T^{-1}\big[U\big] -\mu^{n-1}\times\nu_x\big[U\big]\,\big|\le R(n,L,x).
\]
The proof is now completed by letting $n\to\infty$, then $L\to\infty$ as before.\ $\Box$\\
\\
Here's another proof of Theorem \ref{thmrv}.\\
\\
{\em Second Proof of Theorem \ref{thmrv}:}\, The measure $\mu_{n,x}$ has density with respect to $\mu^n$ given by
\[
f_1(\x)=\frac{\mathbbm{1}\{\sum_{i=1}^n x_i>x\}}{\bar{F}_n(x)}.
\]
The measure
\[
\mu_{n,x}^*=\frac{1}{n}\sum_{j=1}^n\sigma^j\left(\mu^{n-1}\times\nu_x\right)
\]
has density with respect to $\mu^n$ given by
\[
f_2(\x)=\frac{1}{n}\sum_{j=1}^n\frac{\mathbbm{1}\{x_j>x\}}{\bar{F}(x)}=\frac{N_n(\x)}{n\bar{F}(x)},
\]
where $N_n(\x)$ stands for the number of coordinates in $\x$ that are greater than $x$. Now,
\begin{align*}
\|\mu_{n,x}-\mu_{n,x}^*\|_{t.v.}=\int |f_1(\x)-f_2(\x)|\,d\mu(\x)
=\int\left|\frac{\mathbbm{1}\{S_n>x\}}{\bar{F}_n(x)}-\frac{N_n(\X)}{n\bar{F}(x)}\right|d\PP.
\end{align*}
If $n$ is large enough so that
\[
\sup_{x\ge d_n}\left|\frac{\bar{F}_n(x)}{n\bar{F}(x)}-1\right|<\frac{1}{2},
\]
It is a matter of a straightforward computation to see that the preceding expression becomes
\[
\|\mu_{n,x}-\mu_{n,x}^*\|_{t.v.}=2\left(1-\frac{\bar{F}_n(x)}{n\bar{F}(x)}\right)^+\mu_{n,x}\big[N_n(\x)=1\big]+2\,\frac{\PP\big[S_n>x, M_n\le x\big]}{\bar{F}_n(x)}.
\]
The first term on the right hand side above clearly goes to zero, uniformly on $\{x\ge d_n\}$ as $n\to\infty$.
The uniform convergence to zero of the second term  can be deduced from the arguments in the proof of Theorem 4.1 in \cite{MR1652936}, but we give a proof here for the sake of completeness.
\begin{align}
\bar{F}_n(x)&=\PP\big[S_n>x,\, M_n\le x\big]+\PP\big[S_n>x,\, N_n(\X)\ge 1\big] \notag\\
&=\PP\big[S_n>x,\, M_n\le x\big]+n\bar{F}(x)\int_{\sum x_i>x} \frac{d\mu_{n,x}^*(\x)}{N_n(\x)} \notag\\
&= \PP\big[S_n>x,\, M_n\le x\big]+ n\bar{F}(x)\big(1-R(n,x)\big),
\label{max}
\end{align}
where
\begin{align*}
0\le R(n,x)&=\mu_{n,x}^*\big[\sum x_i\le x\big] + \int_{\sum x_i>x} \left( 1-\frac{1}{N_n(\x)}\right)\ d\mu_{n,x}^*(\x)\\
&\le \mu_{n,x}^*\big[\sum x_i\le x\big] +  \mu_{n,x}^*\big[N_n(\x)\ge 2\big]\\
&= \mu_{n,x}^*\big[\sum x_i\le x\big] +\PP\big[M_{n-1}>x\big].
\end{align*}
In view of (\ref{tighty}), it is enough to show that $\mu_{n,x}^*\big[\sum x_i\le x\big]\to 0$, uniformly on $x\ge d_n\vee \ell_n$.
Now,
\begin{align*}
\mu_{n,x}^*\big[\sum x_i\le x\big]&=\int_{S_{n-1}<0}\left(1-\frac{\bar{F}(x-S_{n-1})}{\bar{F}(x)}\right)d\PP\\
&\le D_n(L)+\PP\big[S_{n-1}<-Lb_n\big],
\end{align*}
and this bound goes to zero if we let $n\to\infty$, then $L\to\infty$ by (\ref{ratio}) and (\ref{rc}).  It follows from
(\ref{asymptotic}) and (\ref{max}) that
\[
\lim_{n\to \infty} \sup_{x\ge q_n} \frac{\PP\big[S_n>x, M_n\le x\big]}{\bar{F}_n(x)}=0.
\]\,
\hspace{14.5cm}$\Box$
\section{Some related results}
Let us begin by observing that the assertions of Theorems \ref{thmgen} and  \ref{thmrv} remain valid if we keep $n$ fixed and let $x\to\infty$,
\begin{prop}
If $\mu$ is $\Delta$-subexponential, then for any $n\in\NN$
\begin{equation}
\lim_{x\to\infty}\sup_{A\in\B(\RR^{n-1})}\left|\mu_{n,x}^{\Delta}\circ
T^{-1}\big[A\times\RR \big]-\mu^{n-1}\big[A\big]\right|=0.
\label{sisko}
\end{equation}
\label{thmgenfixed}
\end{prop}
\noindent
We recall the notation $\mu_{n,x}=\mu_{n,x}^{\Delta}$ when $\Delta=(0,\infty)$.
\begin{prop} If $\mu$ is subexponential and $\Delta=(0,\infty)$, then for any $n\in\NN$
\[
\lim_{x\to\infty} \left\|\mu_{n,x}\circ T^{-1}- (\mu^{n-1}\times \nu_x)\right\|_{t.v.}=0.
\]
\label{thmrvfixed}
\end{prop}
\noindent
The proofs of Propositions \ref{thmgenfixed} and \ref{thmrvfixed} are essentially the same as those of Theorems \ref{thmgen} and \ref{thmrv}. We can take $b_n=1$, and instead of using (\ref{localasymptotic}) and (\ref{uct})  we may use the defining relations of ($\Delta$)-subexponentiality. Note that (\ref{sisko}) is proved in \cite{MR2348789} for a family of discrete distributions that includes those with regularly varying tails, subject to $\{S_n=x\}$.\\
\\
So far we have assumed that the interval $\Delta$ is fixed.
We would like to explore what happens when $\Delta=\big(0,s(n,x)\big]$ is a finite interval, but its size $s(n,x)$ grows to infinity with $n$ and $x$. In particular, we would like to understand how the fluctuations of the maximum switch from a stable to a residual life time nature and how fast
 $|\Delta|$ would need to grow for the maximum to be asymptotically independent of the other variables.\\
\\
Denote by $\cal{D}$ the family of semi--open intervals with one
endpoint at 0 contained in the positive half line,
\[
{\cal D}=\Big\{\Delta\subseteq (0,\infty),\, \Delta=(0,z],\,
z>0\Big\}.
\]
Assume that $\mu$ satisfies (\ref{locallongtail}) and
(\ref{localasymptotic}) for some fixed interval $\Delta_0=(0,s_0]$.
Let $\rho_n$ be any sequence increasing to infinity.
\begin{lemma}
If $\mu$ satisfies (\ref{uct}) and (\ref{localasymptotic}) for some
fixed $\Delta_0$, then as  $n\to\infty$ we have
\begin{align*}
\sup_{\Delta \in {\cal D},|\Delta|\ge\rho_n}\,\sup_{x\ge
m_n}\sup_{0\le y\le
b_n}\left|\frac{\mu\big[x-y+\Delta\big]}{\mu\big[x+\Delta\big]}-1\right|\to0 \\
\intertext{and} \sup_{\Delta \in {\cal
D},|\Delta|\ge\rho_n}\,\sup_{x\ge d_n\vee
m_n}\left|\frac{\PP\big[S_n\in
x+\Delta\big]}{n\mu\big[x+\Delta]}-1\right|\to 0.
\end{align*}
\label{lemma}
\end{lemma}
\noindent
{\em Proof:} Let us first note that for any $k\in\NN\cup\{\infty\}$
\begin{equation}
\sup_{x\ge m_n}\sup_{0\le y\le b_n}\left|\frac{\mu\big[x-y+k\Delta_0\big]}{\mu\big[x+k\Delta_0\big]}-1\right|\le \sup_{x\ge m_n}\sup_{0\le y\le b_n}\left|\frac{\mu\big[x-y+\Delta_0\big]}{\mu\big[x+\Delta_0\big]}-1\right|.
\label{easyest}
\end{equation}
In order to see this we split the interval $k\Delta_0=(0,ks_0]$ into $k$ disjoint intervals of length $|\Delta_0|$, $\Delta_0^i=\big((i-1)s_0,i s_0\big]=(i-1)s_0+\Delta_0$, $1\le i\le k$. We get
\begin{align*}
\left|\frac{\mu\big[x-y+k\Delta_0\big]}{\mu\big[x+k\Delta_0\big]}-1\right|&=\frac{1}{\mu\big[x+k\Delta_0\big]}\left|\sum_{i=1}^k \mu\big[x-y+\Delta_0^i]-\mu\big[x+\Delta_0^i\big]\right|\\
&=\frac{1}{\mu\big[x+k\Delta_0\big]}\left|\sum_{i=1}^k \mu\big[x+\Delta_0^i\big] \left( \frac{\mu\big[x+(i-1)s_0-y+\Delta_0\big]}{\mu\big[x+(i-1)s_0+\Delta_0\big]}-1\right)\right|\\
&\le \frac{1}{\mu\big[x+k\Delta_0\big]}\left(\sum_{i=1}^k \mu\big[x+\Delta_0^i\big]\,\, \sup_{x\ge m_n}\sup_{0\le y\le b_n}\left|\frac{\mu\big[x-y+\Delta_0\big]}{\mu\big[x+\Delta_0\big]}-1\right|\,\right)\\
&=\sup_{x\ge m_n}\sup_{0\le y\le b_n}\left|\frac{\mu\big[x-y+\Delta_0\big]}{\mu\big[x+\Delta_0\big]}-1\right|.
\end{align*}
Now, for an arbitrary interval $\Delta\in{\cal D}$ with $|\Delta|\ge \rho_n$, let $k=\big[|\Delta|/|\Delta_0|\big]$ and $\bar{s}=|\Delta|-k|\Delta_0|\le s_0$. If $\bar{\Delta}=(0,\bar{s}]$ we have
\begin{align*}
\frac{\mu\big[x-y+\Delta\big]}{\mu\big[x+\Delta\big]}&=\frac{\mu\big[x-y+\bar{\Delta}\big]+\mu\big[x+\bar{s}-y+k\Delta_0\big]}{\mu\big[x+\bar{\Delta}\big]+\mu\big[x+\bar{s}+k\Delta_0\big]}.
\end{align*}
By dividing both terms of the fraction in the right hand side by $\mu\big[x+\bar{s}+k\Delta_0\big]$, and using (\ref{uct}) and (\ref{easyest}) we see that in order to prove the first assertion of the lemma it suffices to show that
\begin{equation}
\frac{\mu\big[x+\Delta_0\big]}{\mu\big[x+k\Delta_0\big]}\to 0
\label{rest}
\end{equation}
uniformly on $x\ge m_n$, as $k\to\infty$. For this, we may assume without loss of generality that
$k|\Delta_0|\le b_n$ (as otherwise the denominator would be even larger), in which case a computation similar to the previous one
yields
\[
\mu\big[x+k\Delta_0\big]= \sum_{i=1}^k
\mu\big[x+(i-1)s_0+\Delta_0\big]
\ge k\mu\big[x+\Delta_0\big](1-\delta_n),
\]
where $\delta_n$ is the supremum in the right hand side of
(\ref{easyest}), and (\ref{rest}) follows.\\
\\
\noindent For the second assertion we may check again that for any
$k\in\NN\cup\{\infty\}$
\[
\sup_{x\ge d_n}\left|\frac{\PP\big[S_n\in
x+k\Delta_0\big]}{n\mu\big[x+k\Delta_0\big]}-1\right|\le \sup_{x\ge
d_n}\left|\frac{\PP\big[S_n\in
x+\Delta_0\big]}{n\mu\big[x+\Delta_0\big]}-1\right|
\]
and use (\ref{localasymptotic}) and (\ref{rest}) to conclude the proof. $\Box$\\
\\
In view of Lemma \ref{lemma} we can repeat the argument in the proof of Theorem \ref{thmgen} to establish the following.
\begin{prop}
Suppose $\mu$ is $\Delta_0$-subexponential for some finite interval $\Delta_0$. Then,
\[
\lim_{n\to\infty}\sup_{\Delta\in {\cal D}, |\Delta|
\ge\rho_n}\,\sup_{x\ge\,
q_n}\sup_{A\in\B(\RR^{n-1})}\left|\mu_{n,x}^{\Delta}\circ
T^{-1}\big[A\times\RR \big]-\mu^{n-1}\big[A\big]\right|=0,
\]
where the sequence $q_n$ is the same appearing in Theorem 1.
\label{thmgenlargeD}
\end{prop}
\noindent
Theorem \ref{thmrv} also admits a generalisation in this case. Denote by $\nu_x^{\Delta}$ the conditional distribution of $X_i$ subject to $X_i\in x+\Delta$, that is
\[
\nu_x^{\Delta}\big[A\big]=\frac{\mu\big[A\cap(x+\Delta)\big]}{\mu\big[x+\Delta\big]}.
\]
\begin{prop} Suppose $\mu$ is $\Delta_0$-subexponential for some finite interval $\Delta_0$, and let $b_n$ be a sequence such that $S_{n-1}/b_{n}$ is tight. Then
\begin{equation}
\lim_{n\to\infty}\sup_{\Delta \in {\cal D},|\Delta|\ge
\rho_nb_{n}}\sup_{x\ge q_n} \left\|\mu_{n,x}^{\Delta}\circ
T^{-1}- (\mu^{n-1}\times \nu_x^{\Delta})\right\|_{t.v.}=0.
\label{thm2assertlargeD}
\end{equation}
\label{thmrvlargeD}
\end{prop}
\noindent The proof is the same as that of Theorem \ref{thmrv}, if
we substitute $\bar{F}(x),\ \bar{F}_n(x)$ and $I$ by
$\mu\big[x+\Delta\big],\ \PP\big[S_n\in x+\Delta\big]$ and
$\tilde{I}=(x+Lb_n,x+|\Delta|-Lb_n]$ respectively. To show that
$R(n,L,x)\to 0$ we need to use Lemma \ref{lemma}  and the fact that
$\nu_x^\Delta\big[\tilde{I}'\big]\to 0$, which can be proved just as (\ref{rest}). A
direct consequence is that the maximum becomes asymptotically independent
from the rest of the variables as long as $|\Delta|$ grows faster than $b_{n}$.\\
\\
We now discuss the fluctuations of the maximum
$M_n=\max_{1\le i\le n} X_i$ subject to $S_n\in
x+\Delta$, when the distribution $\mu$ is both in
the sum-domain of attraction of a stable law $H$, that is, there
exists a sequence $b_n \uparrow \infty$ such that $S_{n-1}/b_{n}\to
H$, and the max-domain of attraction of an extreme value
distribution (necessarily Gumbel or Fr\'echet). In most cases these can be easily derived from the
preceding results. The following diagram illustrates the typical behaviour for locally subexponential distributions arising in
applications. We assume throughout that $x\ge q_n$, and $\rho_n$ is a sequence growing to infinity.\\
\\
\begin{center}
\begin{picture}(0,0)%
\includegraphics{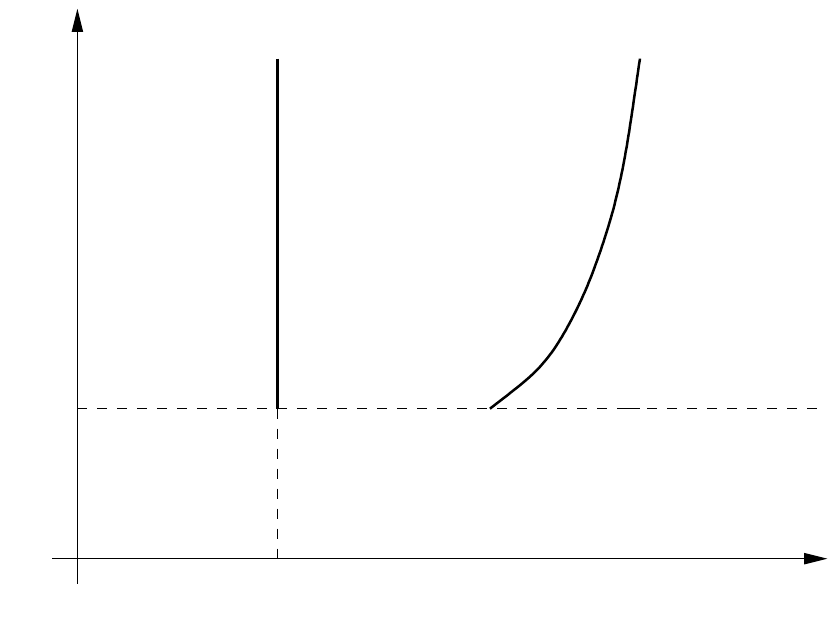}%
\end{picture}%
\setlength{\unitlength}{3158sp}%
\begingroup\makeatletter\ifx\SetFigFont\undefined%
\gdef\SetFigFont#1#2#3#4#5{%
  \reset@font\fontsize{#1}{#2pt}%
  \fontfamily{#3}\fontseries{#4}\fontshape{#5}%
  \selectfont}%
\fi\endgroup%
\begin{picture}(4977,3711)(136,-3289)
\put(2326,-811){\makebox(0,0)[lb]{\smash{{\SetFigFont{8}{9.6}{\rmdefault}{\mddefault}{\updefault}$\frac{M_n-x}{|\Delta|}\to U$}}}}
\put(151,-2011){\makebox(0,0)[lb]{\smash{{\SetFigFont{10}{12.0}{\rmdefault}{\mddefault}{\updefault}$q_n$}}}}
\put(4726,-3211){\makebox(0,0)[lb]{\smash{{\SetFigFont{10}{12.0}{\rmdefault}{\mddefault}{\updefault}$|\Delta|$}}}}
\put(3601,239){\makebox(0,0)[lb]{\smash{{\SetFigFont{10}{12.0}{\rmdefault}{\mddefault}{\updefault}$|\Delta|=\psi(x)$}}}}
\put(226,239){\makebox(0,0)[lb]{\smash{{\SetFigFont{10}{12.0}{\rmdefault}{\mddefault}{\updefault}$x$}}}}
\put(1651,-3211){\makebox(0,0)[lb]{\smash{{\SetFigFont{10}{12.0}{\rmdefault}{\mddefault}{\updefault}$b_{n}$}}}}
\put(3901,-1411){\makebox(0,0)[lb]{\smash{{\SetFigFont{10}{12.0}{\rmdefault}{\mddefault}{\updefault}$\frac{M_n-x}{\psi(x)}\to\Lambda$}}}}
\put(676,-286){\makebox(0,0)[lb]{\smash{{\SetFigFont{8}{9.6}{\rmdefault}{\mddefault}{\updefault}$\frac{M_n-x}{b_{n}}\to -H$}}}}
\end{picture}%

\end{center}

\noindent {\em When} $|\Delta| \le b_{n}/\rho_n$ the maximum
fluctuations are the same as the ones that arise when conditioning
on a finite interval. This holds for all $\Delta_0$-subexponential $\mu$ in the domain of
attraction of a stable law. It follows from Proposition \ref{thmgenlargeD}
and the fact that $\mu_{n,x}^{\Delta}\big[M_n+\sum_{j=1}^{n-1}(TX)_j\in x+\Delta\big]=1$.\\
\\
{\em When} $|\Delta|\ge \rho_n\psi(x)$ the maximum fluctuations are the same as those
obtained by conditioning on $\{S_n>x\}$. This holds for all $\Delta_0$-subexponential $\mu$ in the max-domain of attraction
of the Gumbel or the Fr\'echet distribution. Note that in view of Lemma \ref{lemma}, (\ref{locallongtail}) implies that
\[
\sup_{x\ge m_n}\sup_{|y|\le b_n}\left|\frac{\bar{F}(x-y)}{\bar{F}(x)}-1\right|\longrightarrow 0.
\]
Recall that if the function $\psi$ and the distribution $\Lambda$ are those appearing in (\ref{maxfluct}), for all $t$ for which $\bar{\Lambda}(t)$ is defined we have
\begin{equation}
\frac{\bar{F}\big(x+t\psi(x)\big)}{\bar{F}(x)}\longrightarrow \bar{\Lambda}(t),
\label{Frembel}
\end{equation}
so we must have $\inf_{x\ge m_n}\psi(x)/{b_n}\to\infty.$ Since by Proposition \ref{thmrvlargeD}
\[
\mu_{n,x}^\Delta\Big[\frac{M_n-x}{\psi(x)}\le
t\Big]-\frac{\bar{F}(x)-\bar{F}\big(x+t\psi(x)\big)}{\bar{F}(x)-\bar{F}(x+|\Delta|)}\longrightarrow
0 \text{ as }n\to\infty,
\]
uniformly on $x\ge q_n$, $|\Delta|\gg b_n$, it follows from (\ref{Frembel}) that $(M_n-x)/\psi(x)$ converges to $\Lambda$.\\
\\
{\em In the intermediate region} $b_n\rho_n\le|\Delta|\le\psi(x)/\rho_n$ we prove below that $(M_n-x)/|\Delta|$
converges to a uniform random variable $U$ in $[0,1]$ under some monotonicity condition on $\mu$. Precisely, this holds if $\mu$ is $\Delta_0$-subexponential and in the max-domain of attraction of the Gumbel or Fr\'echet distribution, $S_{n-1}/b_n$ is tight,
and furthermore there exists an interval $\Delta_1=[0,s_1)$ such that $\mu[x+\Delta_1]$ is eventually decreasing in $x$.
We then provide an example to show that the result does not hold in general when the last condition is violated. With this caveat, note that the subexponential distributions typically used in modelling do satisfy all the above conditions.\\
\\
The argument is an adaptation of the proof of Theorem 3.10.11 in \cite{MR1015093}. For any $0<\eps<\half$, if we set $k=\big[\frac{\eps\psi(x)}{|\Delta_1|}\big]$ we have for all sufficiently large $x$
\[
\bar{F}(x)-\bar{F}\big(x+\eps \psi(x)\big)\le \sum_{i=0}^{k} \mu\big[x+is_1+\Delta_1\big]\le (k+1)\mu\big[x+\Delta_1\big],
\]
by monotonicity. Likewise,
\[
\bar{F}\big(x-\eps\psi(x)\big)-\bar{F}\big(x+|\Delta_1|\big)\ge (k+1)\mu\big[x+\Delta_1\big].
\]
Applying these estimates for $\mu\big[x+\Delta_1\big]$ and $\mu\big[x-y+\Delta_1\big]$ for any $y\in\RR$ we have
\begin{equation}
\frac{\bar{F}(x-y)-\bar{F}\big(x-y+\eps\psi(x)\big)}{\bar{F}\big(x-\eps\psi(x)\big)-\bar{F}\big(x+|\Delta_1|\big)}\le \frac{\mu\big[x-y+\Delta_1\big]}{\mu\big[x+\Delta_1\big]}\le \frac{\bar{F}\big(x-y-\eps\psi(x)\big)-\bar{F}\big(x-y+|\Delta_1|\big)}{\bar{F}(x)-\bar{F}\big(x+\eps\psi(x)\big)}.
\label{derivative}
\end{equation}
Since the left hand side in (\ref{Frembel}) is decreasing in $t$, and $\bar{\Lambda}(\cdot)$ is continuous ($\bar{\Lambda}(t)=e^{-t}$ if $\mu$ is attracted to the Gumbel distribution, and $\bar{\Lambda}(t)=(1+t)^{-\alpha}$ if $\mu$ is attracted to the Fr\'echet distribution),  the convergence in (\ref{Frembel}) is uniform in $t\in[-2\eps,2\eps]$. Hence, 
\[
\lim_{n\to\infty}\sup_{x\ge m_n}\sup_{|y|\le\frac{\psi(x)}{\rho_n}}\left|\frac{\mu\big[x-y+\Delta_1\big]}{\mu\big[x+\Delta_1\big]}-1\right|\le \frac{\bar{\Lambda}(\eps)+\bar{\Lambda}(-\eps)-2\bar{\Lambda}(0)}{\bar{\Lambda}(0)-\bar{\Lambda}(\eps)}.
\]
We can let $\eps\to 0$ now to get
\[
\sup_{x\ge m_n}\sup_{|y|\le\frac{\psi(x)}{\rho_n}}\left|\frac{\mu\big[x-y+\Delta_1\big]}{\mu\big[x+\Delta_1\big]}-1\right|=0.
\]
Just as in Lemma \ref{lemma} this in turn implies that
\[
\sup_{|y|\le|\Delta|}\left|\frac{\mu\big[x-y+u\Delta\big]}{\mu\big[x+u\Delta\big]}-1\right|=0,\text{ as } n\to\infty,
\]
for any $u\in[0,1]$, uniformly of $x\ge m_n$, $\rho_n\le|\Delta|\le\psi(x)/\rho_n$, and therefore
\[
\nu_x^{\Delta}\big[x+u\Delta\big]=\frac{\mu\big[x+u\Delta\big]}{\mu\big[x+\Delta\big]}\longrightarrow u.
\]
Now by Proposition \ref{thmrvlargeD} the distribution of the maximum $M_n$ is asymptotically equal to $\nu_x^{\Delta}$ when $|\Delta|\ge \rho_n b_n$, hence $(M_n-x)/|\Delta|\stackrel{d}{\longrightarrow} U$. $\Box$\\
\\
This result is not generally valid without some regularity
assumption on $\mu$. In order to see this, let $\nu$ be the Pareto distribution with density function
\[
\phi(x)=\frac{\alpha}{x^{\alpha+1}}\qquad x\ge 1,
\]
where $\alpha>0$. Consider increasing sequences $c_k\to\infty$ with $c_{k+1}/c_k\to 1$,
$d_k=\sum_{j=1}^k c_{j-1}$, $k\ge 1$, some fixed $0<\epsilon<1$, and
define a distribution $\mu$ with density function
\[
 \bar{\phi}(x)=\sum_{k\ge 0}\frac{1}{\epsilon c_{k}} \,\nu\big[d_k, d_{k+1}\big)\,
 \mathbbm{1}_{[d_k,d_k +\epsilon c_{k})}(x)\,.
\]
In other words, $\mu$ redistributes the mass assigned by $\nu$ to the interval $[d_k, d_{k+1})$ uniformly over the sub--interval
$[d_k, d_k+\epsilon c_{k})$. It is easy to see that $\mu$ satisfies the conditions of Theorem \ref{thmrvlargeD} with
$H$ an $\alpha$--stable law, and that
it belongs to the max--domain of attraction of the Fr\'echet distribution.\\
\\
Let now $k(n)$ be an increasing sequence such that $c_{k(n)}\ge q_n\gg b_n$, and take $\Delta_n=(0,c_{k(n)}]$,
$x_n=d_{k(n)}\gg q_n$. By Theorem \ref{thmrvlargeD}
\[
\lim_{n\to \infty} \mu_{n,x_n}^{\Delta_n}\left[\frac{M_n-x_n}{|\Delta_n|} \in  [\epsilon, 1)\right]=
\lim_{n\to \infty}\frac{\mu\big[x_n+|\Delta_n|\,[\epsilon,1) \big]}{\mu\big[x_n+\Delta_n\big]}=0\neq 1-\epsilon\,,
\]
where $\epsilon$ is the positive value in the definition of $\mu$, and in particular  $M_n-x_n/|\Delta_n|$ does not converge to a uniform distribution. In fact, with little extra effort one can change the measure $\nu$ and choose sequences $x_n$ and
$\Delta_n$ so that $(M_n-x_n)/|\Delta_n|$ converges to any given distribution.\\
\\
One may even determine the behaviour of the fluctuations on the
critical scales of $|\Delta|$.\\
\\
{\em When} $|\Delta|=a\psi(x)$, $(M_n-x)/\psi(x)$ converges to $\Lambda$ conditioned to being less than $a$. This holds under the conditions required for the case $|\Delta|\gg\psi(x)$ and the proof is also the same.\\
\\
{\em Finally, if} $|\Delta|=ab_n$, then $(M_n-x)/b_{n}\stackrel{d}\longrightarrow aU-H$, where the uniform random variable
$U$ is independent of $H$. This holds whenever $\mu$ satisfies (\ref{locallongtail}) and (\ref{localasymptotic}), and $S_{n-1}/b_n$ converges to a stable disribution $H$.\\
\\
To see this suppose that $|\Delta|=ab_n$ and $t\in[0,1]$. Let us begin by observing that
\begin{equation}
\lim_{n\to\infty}\sup_{x\ge m_n}\left|\frac{\mu\big[x+t\Delta\big]}{\mu\big[x+\Delta\big]}-t\right|=0.
\label{arXiv2}
\end{equation}
This follows from (\ref{arXiversion}). Using the fact that $\mu\big[x+t\Delta\big]$ is increasing in $t$ we can show
that convergence is uniform on $t\in[0,1]$.
The rest of the proof follows that  of Theorem \ref{thmgen}.
\begin{align*}
&\PP\big[X_n>x+ub_n,\ S_n\in x+\Delta,\ m_{\X}=n\big]\ge\PP\big[\X^{n-1}\in G,\ X_n>x+ub_n,\ S_n\in x+\Delta\big]\\
&\qquad\qquad=\int_{\X^{n-1}\in G} \mathbbm{1}\big\{ X_n>x+ub_n,\ x-S_{n-1}<X_n\le x-S_{n-1}+ab_n\big\}\ d\PP\\
&\qquad\qquad=\int_{\stackrel{\X^{n-1}\in G}{S_{n-1}<(a-u)b_n}}\Big(\bar{F}\big(x+(ub_n\vee -S_{n-1})\big)-\bar{F}(x+ab_n-S_{n-1})\Big)\ d\PP.
\end{align*}
Now, write the integral above as the sum of the integrals $I_1$ and $I_2$ over the sets $\{S_{n-1}\le -ub_n\}$ and $\{-ub_n< S_{n-1}<(a-u)b_n\}$, respectively. For $x\ge\ell_n$, we can estimate the first of the integrals by
\[
I_1\ge \big(1-D_n(L)\big)\mu\big[x+\Delta\big]\PP\big[\X^{n-1}\in G,\ S_{n-1}\le -ub_n\big],
\]
while the second one can be estimated by
\[
I_2\ge \big(1-D_n(L)\big)\mu\big[x+\Delta\big]\int_{\stackrel{\X^{n-1}\in G}{-ub_n<S_{n-1}<(a-u)b_n}}\frac{\bar{F}(x)-\bar{F}(x+(a-u)b_n-S_{n-1})}{\bar{F}(x)-\bar{F}(x+ab_n)}\ d\PP.
\]
Using the uniform convergence in (\ref{arXiv2}) we may now pass to the limit to get
\begin{align*}
\liminf_{n\to\infty}\inf_{x\ge q_n}\mu_{n,x}^\Delta\Big[\frac{M_n-x}{b_n}>u\Big]&\ge F_H(-u)+\int_{0<\xi+u<a}\left(1-\frac{\xi+u}{a}\right)\ dF_H(\xi)\\
&=\PP\big[aU-H>u\big],
\end{align*}
where $F_H$ above is the distribution function of $H$, $F_H(x)=\PP\big[H\le x\big]$, and $U$ is distributed uniformly on [0,1] and is independent of $H$.
Similarly, we can prove that
\[
\liminf_{n\to\infty}\inf_{x\ge q_n}\mu_{n,x}^\Delta\Big[\frac{M_n-x}{b_n}\le u\Big]\ge\PP\big[aU-H\le u\big],
\]
so $(M_n-x)/b_n$ converges in distribution to $aU-H$. $\Box$\\
\\
{\bf Acknowledgements} The authors would like to thank Amir Dembo,
Stefan Gro{\ss}kinsky, and Richard Samworth for useful comments during the preparation of
this work. ML has been partially supported by the FP7-REGPOT-2009-1
project ``Archimedes Center for Modeling, Analysis and
Computation'', and IA has been partially supported by the
PICT-2008-0315 project ``Probability and Stochastic Processes''.

\small
\bibliographystyle{plain}
\bibliography{subexpALwithlocal}
\end{document}